\documentclass[a4paper,8pt]{article}
\usepackage{amsmath}
\usepackage{amsthm,amsfonts,amssymb}
\usepackage[dvips]{graphicx}
\usepackage[english]{babel}
\usepackage{color}
\usepackage{graphics}
\usepackage{graphpap}

\providecommand{\U}[1]{\protect\rule{.1in}{.1in}}

\providecommand{\U}[1]{\protect\rule{.1in}{.1in}}
\makeatletter
\def\blfootnote{\xdef\@thefnmark{}\@footnotetext}
\makeatother

\begin{document}

\title{Coupon collector's probabilities and generating function: analytic approach}
\author{Andrea Monsellato}
\maketitle

\begin{abstract}
\noindent
Despite the coupon collector's problem has simple probabilistic solution using inclusion/exclusion principle \cite{PolyaUrn}, starting from a particular type of recurrence differential equation
it is used an analytic approach to recover explicit probabilities and generating function.
Note that this kind of recurrence differential equations often appear in Polya Urn context \cite{PolyaUrn}.
\end{abstract}


\section{Introduction}

\noindent
Coupon collector's problem statement is the following: There are $m$ (finite) different coupons to collect.
The coupons are obtained at random, drawing from a box that contains all coupons numbered from $1$ to $m$.
When a coupon is collected the player keeps it, and the drawn coupon is replaced with a new identical coupon.
The target is to collect at least one coupon of every type.\\

\noindent
There are some different approaches to derive explicit probabilities of the problem, i.e. the probabilities to have collected $k\in[1,m]$ coupons
after $n$ draws, for example see \cite{Ross},\cite{Ross2},\cite{PolyaUrn}.\\

\noindent

\section{From master equation to recurrence differential equation}

\noindent
The master equation collector's problems is the following:

\begin{equation}\label{cc:main_eq}
p_{n,k}=p_{n-1,k-1}\frac{m-k+1}{m}+p_{n-1,k}\frac{k}{m}
\end{equation}

\noindent
where $n$ is the step, $k$ is the number of coupons collected and $0<m<+\infty$ is the total number of coupons.
It is obvious that $0\leq k\leq m$.\\

\noindent
Using conditional probabilities, from (\ref{cc:main_eq}), it follows that, for $k\in[1,m]$:

\begin{align*}
E(e^{\lambda X_n})&=E(E(e^{\lambda X_n}|X_{n-1}))=E\left(e^{\lambda (X_{n-1}+1)}\frac{m-X_{n-1}}{m}+e^{\lambda X_{n-1}}\frac{X_{n-1}}{m}\right)=\\&
=E\left(e^{\lambda} e^{\lambda X_{n-1}}\left(1-\frac{X_{n-1}}{m}\right)+e^{\lambda X_{n-1}}\frac{X_{n-1}}{m}\right)
\end{align*}

\noindent
Now set $f_n(e^\lambda)=E(e^{\lambda X_n})$ follows that

\begin{align*}
E(e^{\lambda X_n})&=E\left(e^{\lambda} e^{\lambda X_{n-1}}\left(1-\frac{X_{n-1}}{m}\right)+e^{\lambda X_{n-1}}\frac{X_{n-1}}{m}\right)=\\&
=e^{\lambda}E\left( e^{\lambda X_{n-1}}1-e^{\lambda X_{n-1}}\frac{X_{n-1}}{m}\right)+E\left(e^{\lambda X_{n-1}}\frac{X_{n-1}}{m}\right)=\\&
=e^{\lambda}f_{n-1}(e^\lambda)-\frac{e^{\lambda}}{m}E\left( e^{\lambda X_{n-1}} X_{n-1}\right)+\frac{1}{m}E\left(e^{\lambda X_{n-1}} X_{n-1}\right)=\\&
=e^{\lambda}f_{n-1}(e^\lambda)-\left(\frac{e^{\lambda}-1}{m}\right)E\left( e^{\lambda X_{n-1}} X_{n-1}\right)=\\&
=e^{\lambda}f_{n-1}(e^\lambda)-\left(\frac{e^{\lambda}-1}{m}\right)f'_{n-1}(e^\lambda)
\end{align*}

\noindent
i.e. it obtains that

\begin{align*}
&f_n(e^\lambda)=e^{\lambda}f_{n-1}(e^\lambda)-\left(\frac{e^{\lambda}-1}{m}\right)f'_{n-1}(e^\lambda)
\end{align*}

\noindent
Changing variable $y=e^\lambda$, set $f_n(e^{\lambda})=g_n(y)$, then making derivative w.r.t. $\lambda$

\begin{align}\label{cc:generating_function_eq}
g_n(y)&=y\left[g_{n-1}(y)+\left(\frac{1-y}{m}\right)g'_{n-1}(y)\right]
\end{align}

\noindent
adding the condition $g_0(y)=1$, that equivalent to $f_0(e^\lambda)=1=E(e^{\lambda X_0})$, i.e. at time $0$ no coupon is collected.\\

\noindent
Then following system has to be solved

\begin{align}\label{cc:generating_function_sys}
\left\{
\begin{array}{ll}
g_n(y)=y\left[g_{n-1}(y)+\left(\frac{1-y}{m}\right)g'_{n-1}(y)\right]\qquad \forall n\geq 1\\
g_0(y)=1
\end{array}
\right.
\end{align}

\section{Recovering coupon collector's probabilities}

\noindent
To solve (\ref{cc:generating_function_sys}) the following solution is proposed

\begin{align}\label{cc:generating_function_sol}
g_n(y)=\sum_{k=1}^{n}\frac{y^k}{m^{n-k}}\prod_{h=0}^{k-1}\left(1-\frac{h}{m}\right)a_{n,k}
\end{align}

\noindent
where $a_{n,k}$ coefficients are to be determined.\\

\noindent
Solving recursively the system (\ref{cc:generating_function_sys}) it holds that:\\

\noindent
For $n=1$ it holds that $g_1(y)=ya_{1,1}$ then substituting in (\ref{cc:generating_function_sys}) results that $a_{1,1}=1$.\\

\noindent
For $n=2$, remember that $a_{1,1}=1$, it holds that

\begin{align*}
&\frac{y}{m}a_{2,1}+y^2a_{2,2}\left(1-\frac{1}{m}\right)=\frac{y}{m}a_{1,1}+y^2a_{1,1}\left(1-\frac{1}{m}\right)\\&
\frac{y}{m}a_{2,1}+y^2a_{2,2}\left(1-\frac{1}{m}\right)=\frac{y}{m}+y^2\left(1-\frac{1}{m}\right)\\&
\end{align*}

\noindent
then $a_{2,1}=1$ and $a_{2,2}=1$.\\

\noindent
For $n\geq 3$ the following general treatment is proposed.\\

\noindent
Considering that

\begin{align*}
&g_{n-1}(y)=\sum_{k=1}^{n-1}\frac{y^k}{m^{n-k-1}}\prod_{h=0}^{k-1}\left(1-\frac{h}{m}\right)a_{n-1,k}\\&
g'_{n-1}(y)=\sum_{k=1}^{n-1}\frac{y^{k-1}}{m^{n-k-1}}k\prod_{h=0}^{k-1}\left(1-\frac{h}{m}\right)a_{n-1,k}
\end{align*}

\noindent
then substituting in (\ref{cc:generating_function_sys}) it holds that

\begin{align*}
&\sum_{k=1}^{n}\frac{y^k}{m^{n-k}}\prod_{h=0}^{k-1}\left(1-\frac{h}{m}\right)a_{n,k}=y\left[\sum_{k=1}^{n-1}\frac{y^k}{m^{n-k-1}}\prod_{h=0}^{k-1}\left(1-\frac{h}{m}\right)a_{n-1,k}+
\frac{1-y}{m}\sum_{k=1}^{n-1}\frac{y^{k-1}}{m^{n-k-1}}k\prod_{h=0}^{k-1}\left(1-\frac{h}{m}\right)a_{n-1,k}\right]\\&
\end{align*}

\noindent
Working on the right side of the previous equation

\begin{align*}
\\&\sum_{k=1}^{n-1}\frac{y^{k+1}}{m^{n-k-1}}\prod_{h=0}^{k-1}\left(1-\frac{h}{m}\right)a_{n-1,k}+
y\frac{1-y}{m}\sum_{k=1}^{n-1}\frac{y^{k-1}}{m^{n-k-1}}k\prod_{h=0}^{k-1}\left(1-\frac{h}{m}\right)a_{n-1,k}=\\&
\\&
=\sum_{k=1}^{n-1}\frac{y^{k+1}}{m^{n-k-1}}\prod_{h=0}^{k-1}\left(1-\frac{h}{m}\right)a_{n-1,k}+
\frac{1}{m}\sum_{k=1}^{n-1}\frac{y^{k}k}{m^{n-k-1}}\prod_{h=0}^{k-1}\left(1-\frac{h}{m}\right)a_{n-1,k}
-\frac{1}{m}\sum_{k=1}^{n-1}\frac{y^{k+1}k}{m^{n-k-1}}\prod_{h=0}^{k-1}\left(1-\frac{h}{m}\right)a_{n-1,k}=\\&
\\&
=\sum_{k=1}^{n-1}\frac{y^{k+1}}{m^{n-k-1}}\prod_{h=0}^{k-1}\left(1-\frac{h}{m}\right)\left(1-\frac{k}{m}\right)a_{n-1,k}+
\sum_{k=1}^{n-1}\frac{y^{k}}{m^{n-k}}k\prod_{h=0}^{k-1}\left(1-\frac{h}{m}\right)a_{n-1,k}
\end{align*}

\noindent
It is recovered that

\begin{align}\label{cc:generating_function_sol_explicit}
\sum_{k=1}^{n}\frac{y^k}{m^{n-k}}\prod_{h=0}^{k-1}\left(1-\frac{h}{m}\right)a_{n,k}=\sum_{k=1}^{n-1}\frac{y^{k+1}}{m^{n-k-1}}\prod_{h=0}^{k-1}\left(1-\frac{h}{m}\right)\left(1-\frac{k}{m}\right)a_{n-1,k}+
\sum_{k=1}^{n-1}\frac{y^{k}}{m^{n-k}}k\prod_{h=0}^{k-1}\left(1-\frac{h}{m}\right)a_{n-1,k}
\end{align}

\noindent
Now recovering a general equation for coefficients $a_{n,k}$, from previous equation it holds that

\begin{align*}
&\sum_{k=1}^{n}\frac{y^k}{m^{n-k}}\prod_{h=0}^{k-1}\left(1-\frac{h}{m}\right)a_{n,k}=\sum_{k=1}^{n-1}\frac{y^{k+1}}{m^{n-k-1}}\prod_{h=0}^{k-1}\left(1-\frac{h}{m}\right)\left(1-\frac{k}{m}\right)a_{n-1,k}+
\sum_{k=1}^{n-1}\frac{y^{k}}{m^{n-k}}k\prod_{h=0}^{k-1}\left(1-\frac{h}{m}\right)a_{n-1,k}\\&
\\&
\sum_{k=1}^{n}\frac{y^k}{m^{n-k}}\prod_{h=0}^{k-1}\left(1-\frac{h}{m}\right)a_{n,k}=\sum_{k=1}^{n-1}\frac{y^{k+1}}{m^{n-k-1}}\prod_{h=0}^{k}\left(1-\frac{h}{m}\right)a_{n-1,k}+
\sum_{k=1}^{n-1}\frac{y^{k}}{m^{n-k}}k\prod_{h=0}^{k-1}\left(1-\frac{h}{m}\right)a_{n-1,k}\\&
\\&
\sum_{k=1}^{n}\frac{y^k}{m^{n-k}}\prod_{h=0}^{k-1}\left(1-\frac{h}{m}\right)a_{n,k}=\sum_{k=2}^{n}\frac{y^{k}}{m^{n-k}}\prod_{h=0}^{k-1}\left(1-\frac{h}{m}\right)a_{n-1,k-1}+
\sum_{k=1}^{n-1}\frac{y^{k}}{m^{n-k}}k\prod_{h=0}^{k-1}\left(1-\frac{h}{m}\right)a_{n-1,k}\\&
\end{align*}

\noindent
Splitting the sums, w.r.t. $k$, of both side to apply polynomial identity rule to recover the coefficients $a_{n,k}$, i.e.

\begin{align*}
&\sum_{k=2}^{n-1}\frac{y^k}{m^{n-k}}\prod_{h=0}^{k-1}\left(1-\frac{h}{m}\right)a_{n,k}+\frac{y}{m^{n-1}}\prod_{h=0}^{0}\left(1-\frac{h}{m}\right)a_{n,1}+\frac{y^n}{m^{n-n}}\prod_{h=0}^{n-1}\left(1-\frac{h}{m}\right)a_{n,n}=\\&
=\sum_{k=2}^{n-1}\frac{y^{k}}{m^{n-k}}\prod_{h=0}^{k-1}\left(1-\frac{h}{m}\right)a_{n-1,k-1}+\frac{y^{n}}{m^{n-n}}\prod_{h=0}^{n-1}\left(1-\frac{h}{m}\right)a_{n-1,n-1}+\\&
+\sum_{k=2}^{n-1}\frac{y^{k}}{m^{n-k}}k\prod_{h=0}^{k-1}\left(1-\frac{h}{m}\right)a_{n-1,k}+\frac{y}{m^{n-1}}\prod_{h=0}^{0}\left(1-\frac{h}{m}\right)a_{n-1,1}\\&
\\&
\end{align*}

\noindent
By polynomial equality the following relations hold

\begin{align*}
&a_{n,1}=a_{n-1,1}\\&
\\&
a_{n,n}=a_{n-1,n-1}\\&
\\&
a_{n,k}=a_{n-1,k-1}+ka_{n-1,k}
\end{align*}

\noindent
Now it has obtained that the coefficients $a_{n,k}$ are the Stirling Number of second kind type \cite{aigner}, i.e.

\begin{align*}
a_{n,k}&=\frac{1}{k!}\sum_{j=1}^{k}(-1)^{k-j}\binom{k}{j}j^{n}
\end{align*}

\noindent
Finally it holds that

\begin{align}\label{cc:generating_function_final}
\left\{
\begin{array}{ll}
g_0(y)=1, \qquad n=0\\
g_n(y)=\sum_{k=1}^{n}\frac{y^k}{m^{n-k}}\prod_{h=0}^{k-1}\left(1-\frac{h}{m}\right)\frac{1}{k!}\sum_{j=1}^{k}(-1)^{k-j}\binom{k}{j}j^{n}, \qquad n\leq m\\
g_n(y)=\sum_{k=1}^{m}\frac{y^k}{m^{n-k}}\prod_{h=0}^{k-1}\left(1-\frac{h}{m}\right)\frac{1}{k!}\sum_{j=1}^{k}(-1)^{k-j}\binom{k}{j}j^{n}, \qquad n>m
\end{array}
\right.
\end{align}

\noindent
Also from (\ref{cc:generating_function_final}) using the properties of generating function the explicit form of $p_{n,k}$ is recovered, i.e.

\begin{align}\label{cc:probability_function_final_rude}
\left\{
\begin{array}{ll}
p_{0,k}=1,p_{0,k}=0,\quad \forall k\in[1,m]\\
p_{n,k}=\frac{1}{m^{n-k}}\prod_{h=0}^{k-1}\left(1-\frac{h}{m}\right)\frac{1}{k!}\sum_{j=1}^{k}(-1)^{k-j}\binom{k}{j}j^{n},\qquad \forall n\geq 1,\forall k\in[1,m]\\
\end{array}
\right.
\end{align}

\noindent
if in (\ref{cc:probability_function_final_rude}) the following equality $\frac{1}{m^{n-k}}\prod_{h=0}^{k-1}\left(1-\frac{h}{m}\right)=\frac{1}{m^{n}}\binom{m}{k}k!$ has applied, a simplified form of probabilities $p_{n,k}$ is obtained, i.e.

\begin{align}\label{cc:probability_function_final}
\left\{
\begin{array}{ll}
p_{0,0}=1,p_{0,k}=0,\quad\forall k\in[1,m]\\
p_{n,k}=\frac{\binom{m}{k}}{m^{n}}\sum_{j=1}^{k}(-1)^{k-j}\binom{k}{j}j^{n},\qquad \forall n\geq 1,\forall k\in[1,m]\\
\end{array}
\right.
\end{align}

\newpage
\section{Recovering generating function}

\noindent
Consider the following generating function

\begin{align}\label{cc:generating_function_complete}
G_m(x,y)=\sum_{n=0}^{+\infty}\frac{x^n}{n!}g_n(y)
\end{align}

\noindent
then

\begin{align*}
G_m(x,y)&=\sum_{n=0}^{+\infty}\frac{x^n}{n!}g_n(y)=1+\sum_{n=1}^{+\infty}\frac{x^n}{n!}g_n(y)=\\&
=1+\sum_{n=1}^{m}\frac{x^n}{n!}\sum_{k=1}^{n}\frac{y^k}{m^{n-k}}\prod_{h=0}^{k-1}\left(1-\frac{h}{m}\right)a_{n,k}+
\sum_{n=m+1}^{+\infty}\frac{x^n}{n!}\sum_{k=1}^{m}\frac{y^k}{m^{n-k}}\prod_{h=0}^{k-1}\left(1-\frac{h}{m}\right)a_{n,k}=\\&
=1+\sum_{n=1}^{m}\frac{x^n}{n!}\sum_{k=1}^{n}\frac{y^k}{m^{n}}\binom{m}{k}k!a_{n,k}+
\sum_{n=m+1}^{+\infty}\frac{x^n}{n!}\sum_{k=1}^{m}\frac{y^k}{m^{n}}\binom{m}{k}k!a_{n,k}=\\&
=1+\sum_{n=1}^{m}\frac{x^n}{n!}\sum_{k=1}^{n}\frac{y^k}{m^{n}}\binom{m}{k}k!a_{n,k}+
\sum_{k=1}^{m}y^k\binom{m}{k}k!\sum_{n=m+1}^{+\infty}\frac{x^n}{m^{n}n!}a_{n,k}=\\&
=1+\sum_{n=1}^{m}\frac{x^n}{n!}\sum_{k=1}^{n}\frac{y^k}{m^{n}}\binom{m}{k}k!a_{n,k}+
\sum_{k=1}^{m}y^k\binom{m}{k}k!\left[\sum_{n=1}^{+\infty}\frac{x^n}{m^{n}n!}a_{n,k}-\sum_{n=1}^{m}\frac{x^n}{m^{n}n!}a_{n,k}\right]=\\&
\end{align*}

\noindent
using the properties $a_{n,k}=0$ for $k>n$, it holds that

\begin{align*}
&=1+\sum_{k=1}^{n}y^k\binom{m}{k}k!\sum_{n=1}^{m}\frac{x^n}{m^n n!}a_{n,k}+
\sum_{n=1}^{m}\frac{x^n}{n!}\sum_{k=1}^{n}\frac{y^k}{m^{n}}\binom{m}{k}k!a_{n,k}-\sum_{n=1}^{m}\frac{x^n}{m^{n}n!}\sum_{k=1}^{m}y^k\binom{m}{k}k!a_{n,k}=\\&
=1+\sum_{k=1}^{m}y^k\binom{m}{k}k!\sum_{n=1}^{+\infty}\frac{x^n}{m^n n!}a_{n,k}
-\sum_{n=1}^{m}\frac{x^n}{m^{n}n!}\sum_{k=n+1}^{m}y^k\binom{m}{k}k!a_{n,k}=\\&
=1+\sum_{k=1}^{m}y^k\binom{m}{k}k!\sum_{n=1}^{+\infty}\frac{x^n}{m^n n!}a_{n,k}
\end{align*}

\noindent
now substituting the $a_{n,k}$ expression, a closed expression of generating function is retrieved

\begin{align}\label{cc:generating_function_complete_final}
\notag G_m(x,y)&=1+\sum_{k=1}^{m}y^k\binom{m}{k}k!\sum_{n=1}^{+\infty}\frac{x^n}{m^n n!}a_{n,k}=\\&
\notag =1+\sum_{k=1}^{m}y^k\binom{m}{k}k!\sum_{n=1}^{+\infty}\frac{x^n}{m^n n!}\frac{1}{k!}\sum_{j=1}^{k}(-1)^{k-j}\binom{k}{j}j^{n}=\\&
\notag =1+\sum_{k=1}^{m}(-1)^{k}y^k\binom{m}{k}\sum_{j=1}^{k}(-1)^{j}\binom{k}{j}\sum_{n=1}^{+\infty}\frac{x^n}{m^n n!}j^{n}=\\&
\notag =1+\sum_{k=1}^{m}(-1)^{k}y^k\binom{m}{k}\sum_{j=1}^{k}(-1)^{j}\binom{k}{j}\left[e^{\frac{x}{m}j}-1\right]=\\&
\notag =1+\sum_{k=1}^{m}(-1)^{k}y^k\binom{m}{k}\left[\left[\left(1-e^{\frac{x}{m}}\right)^k-1\right]-\left[(1-1)^k-1\right]\right]=\\&
\notag =1+\sum_{k=1}^{m}(-1)^{k}y^k\binom{m}{k}\left(1-e^{\frac{x}{m}}\right)^k=\\&
=\left[1-y\left(1-e^{\frac{x}{m}}\right)\right]^m
\end{align}

\end{document}